\documentclass[a4paper, twoside, reqno]{amsart} 

\usepackage{hyperref}					 
\usepackage[latin1]{inputenc}  
\usepackage{float}
\usepackage{graphicx}

\theoremstyle{definition}

\theoremstyle{remark}

\numberwithin{equation}{section}

\DeclareMathOperator{\sn}{sn}

\DeclareMathOperator{\hav}{hav}
\DeclareMathOperator{\am}{am}

\begin{document}

\title{The Liouville parametrization of a triaxial ellipsoid}
\author{\sc C{\v a}lin--{\c S}erban B{\v a}rbat}
\address{Dieselstr. 19\\80993 München\\Germany} 
\email{calin.barbat@web.de}
\keywords{Liouville surface, ellipsoid, elliptic integral of third kind}
\date{\today}
\begin{abstract}
In this article we will construct the Liouville parametrization of the triaxial ellipsoid. In the literature quadrics are given as examples of Liouvillesurfaces, yet no one gives such a parametrization.
\end{abstract}
\maketitle

\section{Introduction}

\begin{figure}[H]
	\centering
	$\begin{array}{cccc}
   \includegraphics[scale=0.5]{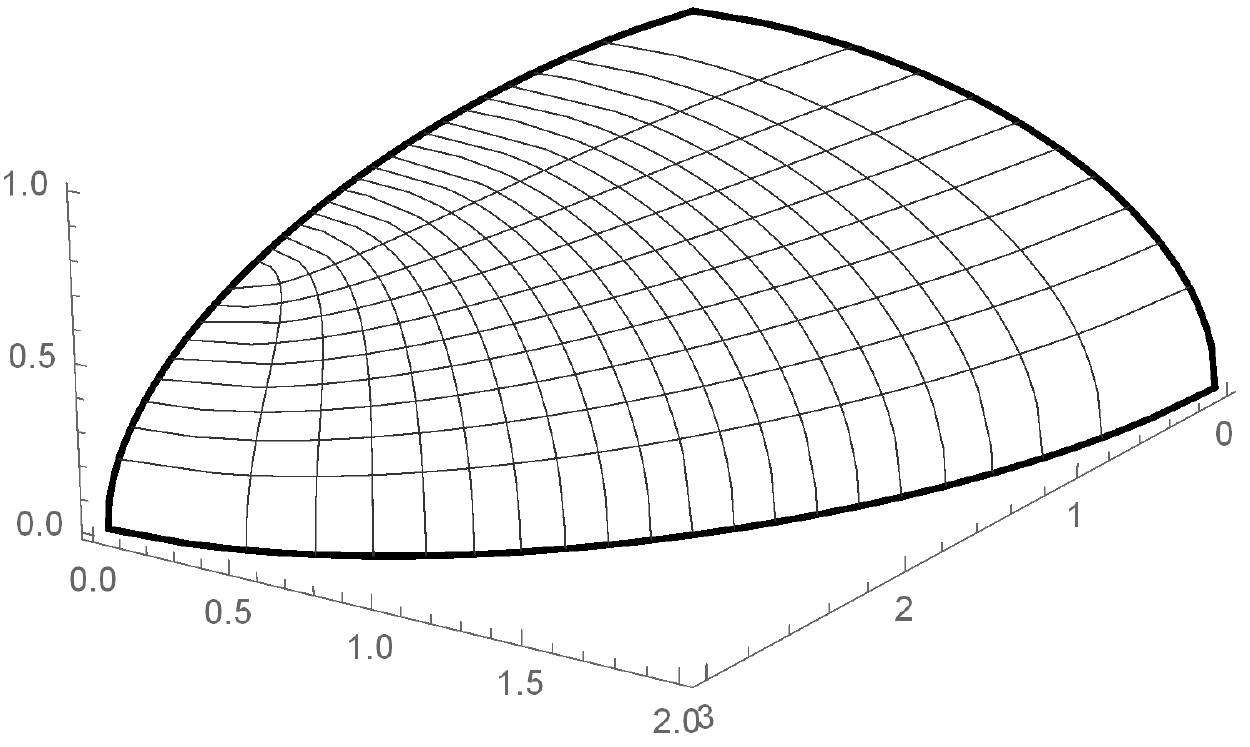} &
   \includegraphics[scale=0.5]{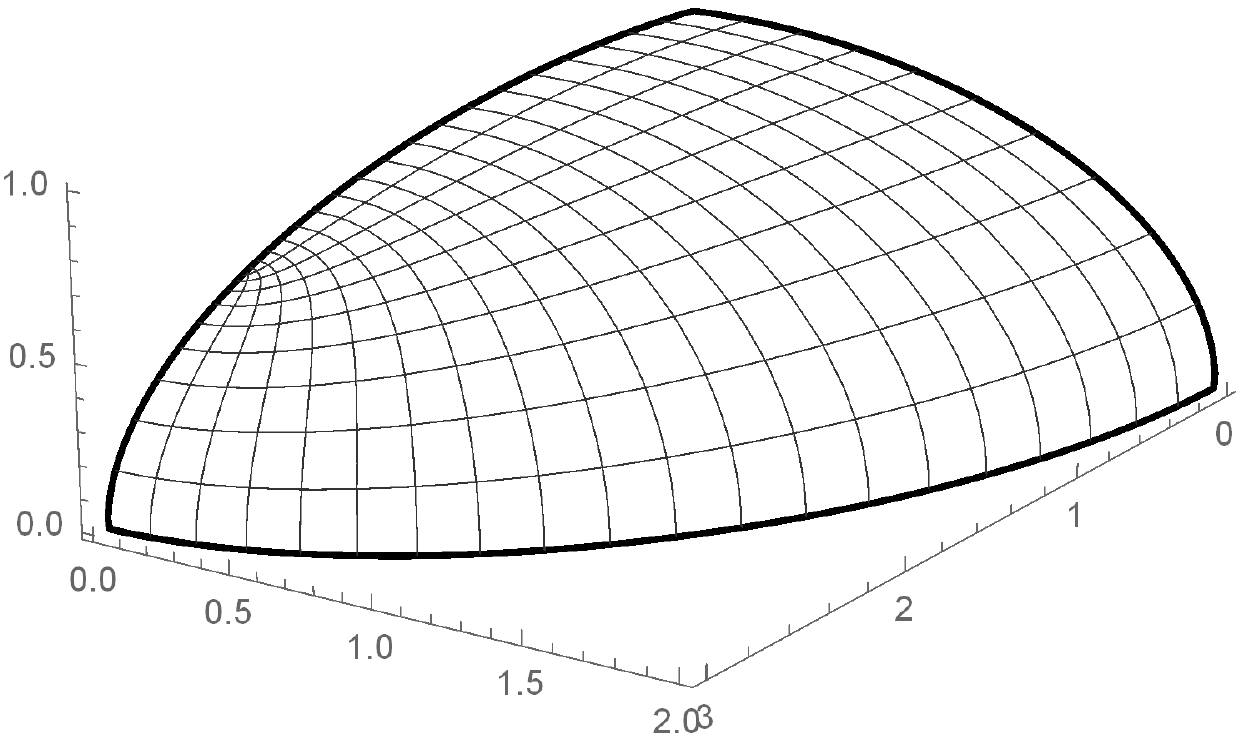}
  \end{array}$
	\caption{Standard curvature line (left image) and Liouville (right image) parametrization of a triaxial ellipsoid}
	\label{fig:result}
\end{figure}

In the literature that you will find at the end of this article (see \cite{bmu}, \cite{nyr}, \cite{val}), the authors describe how to map a triaxial ellipsoid conformally to a plane. The best paper (to my knowledge) on this matter is \cite{nyr} because it actually computes (making use of elliptic integrals) the integrals already given by Jacobi in his ``Lectures on Dynamics''. In this article we want to go in the opposite direction and map a plane rectangle conformally to a triaxial ellipsoid in such a way that the map has a Liouville line element. The result can be seen in the right image of figure \ref{fig:result}.

\section{Standard curvature line parametrization of the triaxial ellipsoid}

We will start here with the standard curvature line parametrization of the triaxial ellipsoid with semi-axes $0 < c < b < a$:
\begin{gather*} 
\text{Ellipsoid}(u,v)=\left( 
\sqrt{\frac{a^2 (a^2 - u) (a^2 - v)}{(a^2 - b^2) (a^2 - c^2)}}, 
\sqrt{\frac{b^2 (b^2 - u) (b^2 - v)}{(b^2 - c^2) (b^2 - a^2)}}, 
\sqrt{\frac{c^2 (c^2 - u) (c^2 - v)}{(c^2 - a^2) (c^2 - b^2)}}
\right)^t 
\end{gather*}
where $0 < c^2 < v < b^2 < u < a^2$ (see left image of figure \ref{fig:result}).

The coefficients of the first fundamental form are computed as follows:
\begin{align*} 
g_{11}(u, v)&=\left\langle \partial_u \text{Ellipsoid}(u,v), \partial_u \text{Ellipsoid}(u,v) \right\rangle=\frac{1}{4} (u-v) f(u)\\
g_{12}(u, v)&=\left\langle \partial_u \text{Ellipsoid}(u,v), \partial_v \text{Ellipsoid}(u,v) \right\rangle=g_{21}(u, v) = 0\\
g_{22}(u, v)&=\left\langle \partial_v \text{Ellipsoid}(u,v), \partial_v \text{Ellipsoid}(u,v) \right\rangle=\frac{1}{4} (u-v) (-f(v))
\end{align*}
with the function $f$ defined as:
\begin{align*} 
f(t)=\frac{t}{(a^2 - t) (b^2 - t) (c^2 - t)}
\end{align*}

The line element of the ellipsoid is:
\begin{align} 
{ds}^2 = g_{11}(u, v) {du}^2 + g_{22}(u, v) {dv}^2 = \frac{1}{4} (u-v) (f(u){du}^2-f(v){dv}^2) \label{old_ds}
\end{align}

\section{Conformal map from ellipsoid to plane}

What we want to achieve is the following Liouville form of this line element \ref{old_ds}:
\begin{align} 
{ds}^2 &= \frac{1}{4} (U(x)-V(y)) ({dx}^2+{dy}^2) \label{new_ds}
\end{align}

If formulas \ref{old_ds} and \ref{new_ds} are to be the same we must have:
\begin{align*} 
       dx = \sqrt{+f(u)}{du} \quad \text{and} \quad
       dy = \sqrt{-f(v)}{dv}
\end{align*}

By integrating we get formulas corresponding to (7) and (8) from \cite{nyr}:
\begin{align*} 
       X(u) &= \int_{b^2}^{u} \sqrt{+f(t)}{dt} = F_1(u)-F_1(b^2) = F_1(u) \\
       Y(v) &= \int_{c^2}^{v} \sqrt{-f(t)}{dt} = F_2(v)-F_2(c^2) = F_2(v)
\end{align*}
with
\begin{align*}
F_1\left(t\right)&=\frac{2 b^2 i}{c \sqrt{ a^2 - b^2 }} \Pi\left(n_1; \varphi_1(t) | m_1 \right)  \\ 
F_2\left(t\right)&=\frac{2 c^2  }{b \sqrt{ a^2 - c^2 }} \Pi\left(n_2; \varphi_2(t) | m_2 \right) 
\end{align*}
where $i=\sqrt{-1}$ and
\begin{align*}
n_1 &= 1 - \frac{b^2}{c^2} &
\varphi_1(t) &= \arcsin \left({-i c \sqrt{\frac{ t - b^2 }{\left( b^2 - c^2 \right) t  }}}\right) 
          & m_1 &= \frac{{a}^2\left( c^2 - b^2 \right) }{{c}^2\left( a^2 - b^2 \right) }\\
n_2 &= 1 - \frac{c^2}{b^2} &
\varphi_2(t) &= \arcsin \left({b \sqrt{\frac{ t - c^2 }{\left( b^2 - c^2 \right) t  }}}\right) 
          & m_2 &= \frac{{a}^2\left( b^2 - c^2 \right) }{{b}^2\left( a^2 - c^2 \right) }
\end{align*}
and the incomplete elliptic integral of the third kind is defined as follows:
\begin{align*}
\Pi(n; \varphi | m)=\int_0^{\varphi} \frac{d\theta}{(1 - n\sin^2 \theta)\sqrt{1 - m \sin^2 \theta}}
\end{align*}

\section{Liouville map from plane to ellipsoid}

We are interested in the functions $U(x)$ and $V(y)$. But we have $X(u)$ and $Y(v)$, which cannot be inverted easily. We have two alternatives:
\begin{enumerate}
\item The first alternative is to define a generalized Jacobi amplitude $\am(n;z|m)$ as inverse function of the elliptic integral of the third kind. That means 
\begin{align*}
z &= \Pi(n; \varphi | m) \\
\am(n; z | m) &= \varphi
\end{align*}
The Jacobi amplitude as special case can be expressed in terms of this generalized Jacobi amplitude as $\am(z|m)=\am(0;z|m)$. With the generalized Jacobi amplitude we can invert the elliptic integrals of the third kind and get:
\begin{align*}
U(x)&=\frac{b^2}{1-n_1 \sin ^2\left(\am \left(n_1;\frac{x c \sqrt{a^2-b^2}}{2 i b^2}|m_1\right)\right)} \\
V(y)&=\frac{c^2}{1-n_2 \sin ^2\left(\am \left(n_2;\frac{y b \sqrt{a^2-c^2}}{2 c^2}  |m_2\right)\right)}
\end{align*}
We can introduce the generalized Jacobi elliptic function $\sn(n;z|m) = \sin(\am(n;z|m))$ and then we have:
\begin{align*}
U(x)&=\frac{b^2}{1-n_1 \sn ^2\left(n_1;\frac{x c \sqrt{a^2-b^2}}{2 i b^2}|m_1\right)} \\
V(y)&=\frac{c^2}{1-n_2 \sn ^2\left(n_2;\frac{y b \sqrt{a^2-c^2}}{2 c^2}  |m_2\right)}
\end{align*}

\item The other alternative is to use a series representation for $X(u)$ and $Y(v)$ and then compute the reverse/inverse series, giving a series representation for $U(x)$ and $V(y)$.
We first expand $X(u)$ in a series about the point $u_0=b^2$ and $Y(v)$ in a series about the point $v_0=c^2$:
\begin{align*}
X(u) &= \sum_{k=0}^{\infty} A_{2k+1} \left(\frac{\sqrt{u-b^2}}{\sqrt{\left(a^2-b^2\right) \left(b^2-c^2\right)}}\right)^{2k+1} \\
Y(v) &= \sum_{k=0}^{\infty} B_{2k+1} \left(\frac{\sqrt{v-c^2}}{\sqrt{\left(a^2-c^2\right) \left(b^2-c^2\right)}}\right)^{2k+1}
\end{align*}
The first three coefficients are:
\begin{align*}
A_1&= 2 b, 																& B_1 &= 2 c \\ 
A_3&=\frac{b^4 - a^2 c^2}{3 b}, 					& B_3 &= -\frac{c^4 - a^2 b^2}{3 c} \\
A_5&=\frac{-a^4 c^4+4 a^4 b^2 c^2-10 a^2 b^4 c^2+4 a^2 b^2 c^4+3 b^8}{20 b^3} \\
B_5&=\frac{-a^4 b^4+4 a^4 b^2 c^2+4 a^2 b^4 c^2-10 a^2 b^2 c^4+3 c^8}{20 c^3}
\end{align*}
By computing the reverse/inverse series we get:
\begin{align*}
U(x) &= b^2 + \left( a^2 - b^2 \right) \left( b^2 - c^2 \right) \sum_{k=1}^{\infty} C_{2k} x^{2k}, \\
V(y) &= c^2 + \left( c^2 - a^2 \right) \left( c^2 - b^2 \right) \sum_{k=1}^{\infty} D_{2k} y^{2k}
\end{align*}
where the first three coefficients are:
\begin{align*}
C_2 &= \frac{1}{4 b^2}, 																  & D_2 &= \frac{1}{4 c^2} \\
C_4 &= -\frac{b^4 - a^2 c^2}{48 b^6},	                    & D_4 &= \frac{ c^4 - a^2 b^2 }{48 c^6} \\
C_6 &= \frac{11 a^4 c^4-9 a^4 b^2 c^2-9 a^2 b^2 c^4+5 a^2 b^4 c^2+2 b^8}{2880 b^{10}},\\
D_6 &= \frac{11 a^4 b^4-9 a^4 b^2 c^2-9 a^2 b^4 c^2+5 a^2 b^2 c^4+2 c^8}{2880 c^{10}}
\end{align*}
\end{enumerate}

Then the Liouville parametrization of the ellipsoid is given by:
\begin{gather*} 
\text{Ellipsoid}(U(x),V(y))
\end{gather*}
where $0 = X(b^2) < x < X(a^2)$ and $0 = Y(c^2) < y < Y(b^2)$.

\section{Differential equations}

If we plug $u=U(x)$ and $v=V(y)$ in the equation \ref{old_ds} of the line element of the ellipsoid we get:
\begin{align*}
{ds}^2 = \frac{1}{4} \left(U(x)-V(y)\right) 
                                \left(f(U(x))\left(\frac{dU(x)}{dx}\right)^2 dx^2
                                -f(V(y))\left(\frac{dV(y)}{dy}\right)^2 dy^2 \right)
\end{align*}
Comparing this formula with \ref{new_ds} we see that the functions $U(x)$ and $V(y)$ satisfy the following differential equations:
\begin{align*}
    f(U(x))\left(\frac{dU(x)}{dx}\right)^2 = +1 \quad \text{and} \quad
    f(V(y))\left(\frac{dV(y)}{dy}\right)^2 = -1
\end{align*}

\section{Remark about the figure}

Because we don't have the complete series $U(x)$ and $V(y)$ (but only an approximation, with a few terms) we need another method for drawing a quite good figure of the Liouville ellipsoid for $a=3$, $b=2$, $c=1$:
\begin{enumerate}
	\item Collect the points $(X(u_k), u_k)$ at the values $b^2 = u_0 < u_1 < \ldots < u_k = \frac{(n-k) u_0 + k u_n}{n} < \ldots < u_n = a^2$.
	\item Interpolate these points with some smooth function (possibly piecewise defined), giving a good approximation $\widetilde{U}(x) \approx U(x)$.
	\item Do the same with the points $(Y(v_k), v_k)$ and get $\widetilde{V}(y) \approx V(y)$.
	\item Draw the figure with the parametrization $\text{Ellipsoid}(\widetilde{U}(x),\widetilde{V}(y))$ for $0 = X(b^2) < x < X(a^2)$ and $0 = Y(c^2) < y < Y(b^2)$.
\end{enumerate}
One could also try to numerically invert the functions $X(u)$ and $Y(v)$ as described in \cite{fuk} to get approximations to $U(x)$ and $V(y)$.

\section{Open problem}

Look at the series representations of the (haversed sine)/haversine and inverse haversine functions:
\begin{align*}
\hav(z)      &= \sin^2\left(\frac{z}{2}\right) = \sum_{k=1}^{\infty} \frac{(-1)^{k-1}}{2 (2 k)!}z^{2 k} \\
             &= \frac{z^2}{4}-\frac{z^4}{48}+\frac{z^6}{1440}-\frac{z^8}{80640}+\frac{z^{10}}{7257600}+\ldots \\
\hav^{-1}(z) &= 2\arcsin\left(\sqrt{z}\right)  = \sum_{k=0}^{\infty} \frac{\binom{2 k}{k}}{2^{2 k-1}(2 k+1)}(\sqrt{z})^{2 k+1} \\
             &= 2 \sqrt{z}+\frac{z^{3/2}}{3}+\frac{3 z^{5/2}}{20}+\frac{5 z^{7/2}}{56}+\frac{35 z^{9/2}}{576}+\ldots
\end{align*}
These coefficients/numbers also appear in the series of $U(x)$, $V(y)$ and $X(u)$, $Y(v)$. The series for $U(x)$ and $V(y)$ can be written now as:
\begin{align*}
U(x) &= b^2 + \left( a^2 - b^2 \right) \left( b^2 - c^2 \right) \sum_{k=1}^{\infty} \frac{(-1)^{k-1}\gamma_{2k}}{2 (2 k)! b^{2(2 k-1)}} x^{2k}, \\
V(y) &= c^2 + \left( c^2 - a^2 \right) \left( c^2 - b^2 \right) \sum_{k=1}^{\infty} \frac{\delta_{2k}}{2 (2 k)! c^{2(2 k-1)}} y^{2k}
\end{align*}
where the first three coefficients are:
\begin{align*}
\gamma_2 &= 1, 													  			  & \delta_2 &= 1 \\
\gamma_4 &= b^4 - a^2 c^2,  	                    & \delta_4 &= c^4 - a^2 b^2 \\
\gamma_6 &= \frac{11 a^4 c^4-9 a^4 b^2 c^2-9 a^2 b^2 c^4+5 a^2 b^4 c^2+2 b^8}{2},\\
\delta_6 &= \frac{11 a^4 b^4-9 a^4 b^2 c^2-9 a^2 b^4 c^2+5 a^2 b^2 c^4+2 c^8}{2}
\end{align*}
The series representations for $X(u)$ and $Y(v)$ are as follows:
\begin{align*}
X(u) &= \sum_{k=0}^{\infty} \frac{\binom{2 k}{k} \alpha_{2k+1}}{2^{2 k-1}(2 k+1)b^{2 k-1}} \left(\frac{\sqrt{u-b^2}}{\sqrt{\left(a^2-b^2\right) \left(b^2-c^2\right)}}\right)^{2k+1} \\
Y(v) &= \sum_{k=0}^{\infty} \frac{\binom{2 k}{k} (-1)^k \beta_{2k+1}}{2^{2 k-1}(2 k+1)c^{2 k-1}}
\left(\frac{\sqrt{v-c^2}}{\sqrt{\left(a^2-c^2\right) \left(b^2-c^2\right)}}\right)^{2k+1}
\end{align*}
The first three coefficients are:
\begin{align*}
\alpha_1&= 1, 																& \beta_1 &= 1\\ 
\alpha_3&=b^4 - a^2 c^2, 					            & \beta_3 &= c^4 - a^2 b^2 \\
\alpha_5&=\frac{-a^4 c^4+4 a^4 b^2 c^2-10 a^2 b^4 c^2+4 a^2 b^2 c^4+3 b^8}{3} \\
\beta_5 &=\frac{-a^4 b^4+4 a^4 b^2 c^2+4 a^2 b^4 c^2-10 a^2 b^2 c^4+3 c^8}{3}
\end{align*}
It is possible to calculate the first coefficients $\alpha_{2k+1}$, $\beta_{2k+1}$, $\gamma_{2k}$ and $\delta_{2k}$ of the series expansions of $X(u)$, $Y(v)$, $U(x)$ and $V(y)$. But I have not been able to get the closed general form of these coefficients. This is an open problem and I would like to hear from you, if you make progress on it.

\section{Acknowledgements}

I want to thank Prof. Maxim Nyrtsov for sending me his paper about the Jacobi conformal map from ellipsoid to plane. I want to thank Albert D. Rich for his invaluable help in computing the two integrals $F_1(t)$ and $F_2(t)$. He will add these integrals to his rule based integrator (see \cite{adr}).
My thanks also go to my family, my friends and my employer for their support.

\end{document}